\documentclass{amsart}
\usepackage{amsfonts}

\numberwithin{equation}{section}

\begin{document}
\title[]{Representation of multivariate functions via the potential theory}
\author[F.-C. C\^ {\i}rstea]{Florica-Corina C\^ {\i}rstea}
\address{School of Computer Science and Mathematics, Victoria University of
Technology, PO Box 14428, Melbourne City MC, Victoria 8001, Australia}
\email{florica@sci.vu.edu.au}
\author[S. S. Dragomir]{Sever Silvestru Dragomir}
\address{School of Computer Science and Mathematics, Victoria University of
Technology, PO Box 14428, Melbourne City MC, Victoria 8001, Australia}
\email{Sever.Dragomir@vu.edu.au}

\begin{abstract}
In this paper, by the use of Potential Theory, some representation results
for multivariate functions from the Sobolev spaces $W^{1,p}(\Omega ),$ in
terms of the double layer potential and the fundamental solution of
Laplace's equation are pointed out. Applications for multivariate
inequalities of Ostrowski type are also provided.
\end{abstract}

\maketitle
\newtheorem{th1}{Theorem} \newtheorem{lem}{Lemma} \newtheorem{prop}{%
Proposition} \newtheorem{cor}{Corollary} \theoremstyle{remark} \newtheorem{%
rem}{Remark}

\section{Introduction}

The following representation for an absolutely continuous function $f:\left[
a,b\right] \rightarrow \mathbb{R}$ in terms of the integral mean is known in
the literature as Montgomery identity 
\begin{equation*}
f\left( x\right) =\frac{1}{b-a}\int_{a}^{b}f\left( t\right) dt+\frac{1}{b-a}%
\int_{a}^{b}p\left( t,x\right) f^{\prime }\left( t\right) dt,\text{ }x\in %
\left[ a,b\right] ;
\end{equation*}
where $p:\left[ a,b\right] ^{2}\rightarrow \mathbb{R,}$ is given by 
\begin{equation}
p\left( t,x\right) =\left\{ 
\begin{array}{cc}
t-a & \text{if }a\leq t\leq x \\ 
t-b & \text{if }x<t\leq b
\end{array}
\right. .  \label{eq1}
\end{equation}
In the last decade, many authors (see for example \cite{DR} and the
references therein) have extended the above result for different classes of
functions defined on a compact interval, including: functions of bounded
variation, monotonic functions, convex functions, $n$-time differentiable
functions whose derivatives are absolutely continuous or satisfy different
convexity properties etc...and pointed out sharp inequalities for the
absolute value of the difference 
\begin{equation*}
D(f;x):=f\left( x\right) -\frac{1}{b-a}\int_{a}^{b}f\left( t\right) dt,\text{
}x\in \left[ a,b\right] .
\end{equation*}
The obtained results have been applied in Approximation Theory, Numerical
Integration, Information Theory and other related domains.

We have, see for instance \cite[p. 2]{DR}, the following \textit{Ostrowski
type inequalities} 
\begin{eqnarray*}
&&\left| D(f;x)\right| \\
&\leq &\left\{ 
\begin{array}{cc}
\left[ \frac{1}{4}+\left( \frac{x-\frac{a+b}{2}}{b-a}\right) ^{2}\right]
\left( b-a\right) \left\| f^{\prime }\right\| _{\infty } & \text{ if }%
f^{\prime }\in L_{\infty }\left[ a,b\right] ; \\ 
\frac{1}{\left( p+1\right) ^{1/p}}\left[ \left( \frac{x-a}{b-a}\right)
^{p+1}+\left( \frac{b-x}{b-a}\right) ^{p+1}\right] ^{1/p}\left( b-a\right)
^{1/p}\left\| f^{\prime }\right\| _{q} & 
\begin{array}{c}
\text{if }f^{\prime }\in L_{q}\left[ a,b\right] \\ 
q>1,\frac{1}{p}+\frac{1}{q}=1
\end{array}
; \\ 
\left[ \frac{1}{2}+\left| \frac{x-\frac{a+b}{2}}{b-a}\right| \right] \left\|
f^{\prime }\right\| _{1}; & 
\end{array}
\right.
\end{eqnarray*}
provided $f$ is absolutely continuous and $L_{r}\left[ a,b\right] \left(
1\leq r\leq \infty \right) $ are the usual Lebesgue spaces. The constants $%
\frac{1}{4},\frac{1}{\left( p+1\right) ^{1/p}}$ and $\frac{1}{2}$ are best
possible in the sense that they cannot be replaced by smaller constants.

If the functions $f:\left[ a,b\right] \times \left[ c,d\right] \rightarrow 
\mathbb{R}$ has the partial derivatives $\frac{\partial f\left( t,s\right) }{%
\partial t},\frac{\partial f\left( t,s\right) }{\partial s},$ and $\frac{%
\partial ^{2}f\left( t,s\right) }{\partial t\partial s}$ continuous on $%
\left[ a,b\right] \times \left[ c,d\right] ,$ then one has the
representation \cite[p. 307]{DR} 
\begin{eqnarray*}
f\left( x,y\right)  &=&\frac{1}{(b-a)(d-c)}\int_{a}^{b}\int_{c}^{d}f\left(
t,s\right) dtds \\
&&+\frac{1}{(b-a)(d-c)}\int_{a}^{b}\int_{c}^{d}p\left( t,x\right) \frac{%
\partial f\left( t,s\right) }{\partial t}dtds \\
&&+\frac{1}{(b-a)(d-c)}\int_{a}^{b}\int_{c}^{d}q\left( s,y\right) \frac{%
\partial f\left( t,s\right) }{\partial s}dtds \\
&&+\frac{1}{(b-a)(d-c)}\int_{a}^{b}\int_{c}^{d}p\left( t,x\right) q\left(
s,y\right) \frac{\partial ^{2}f\left( t,s\right) }{\partial t\partial s}dtds,
\end{eqnarray*}
for each $\left( x,y\right) \in \left[ a,b\right] \times \left[ c,d\right] ,$
where $p$ is defined by (\ref{eq1}) and $q$ is the corresponding kernel for
the interval $\left[ c,d\right] .$

Another representation for $f:\left[ a,b\right] \times \left[ c,d\right]
\rightarrow \mathbb{R}$ is \cite[p. 294]{DR} 
\begin{eqnarray*}
f\left( x,y\right)  &=&\frac{1}{b-a}\int_{a}^{b}f\left( t,y\right) dt+\frac{1%
}{d-c}\int_{c}^{d}f\left( x,s\right) ds \\
&&-\frac{1}{(b-a)(d-c)}\int_{a}^{b}\int_{c}^{d}f\left( t,s\right) dtds \\
&&+\frac{1}{(b-a)(d-c)}\int_{a}^{b}\int_{c}^{d}p\left( t,x\right) q\left(
s,y\right) \frac{\partial ^{2}f\left( t,s\right) }{\partial t\partial s}dtds,
\end{eqnarray*}
for each $\left( x,y\right) \in \left[ a,b\right] \times \left[ c,d\right] ,$
provided $\frac{\partial ^{2}f\left( t,s\right) }{\partial t\partial s}$ is
continuous in $\left[ a,b\right] \times \left[ c,d\right] .$

Different Ostrowski type inequalities for multivariate functions may be
stated, see Chapters 5 \& 6 of \cite{DR}.

In this paper, by the use of Potential Theory, some representation results
for multivariate functions from the Sobolev spaces $W^{1,p}(\Omega ),$ where 
$\Omega $ is an open bounded set with smooth boundary in $\mathbb{R}%
^{N},N\geq 2,p\in (N,\infty ],$ in terms of the double layer potential and
the fundamental solution of Laplace's equation are pointed out. Applications
for multivariate inequalities of Ostrowski type are also provided.

\section{Preliminaries}

For $\Omega\subset {\mathbb R}^N$, we denote by $\overline{\Omega}$ its closure and
by $\partial \Omega$ the boundary of $\Omega$.

By a \emph{vector field} we understand an ${\mathbb R}^N$-valued function on a subset
of ${\mathbb R}^N$. If $Z=(z_1, z_2,\ldots, z_N)$ is a differentiable vector field on
an open set $\Omega\subset {\mathbb R}^N$, the \emph{divergence} of $Z$ on $\Omega$
is defined by 
\begin{equation*}
\mathrm{div}Z=\sum_{i=1}^N \frac{\partial z_i}{\partial x_i}.
\end{equation*}

\begin{prop}[The Divergence Theorem] \label{div}
Let $\Omega\subset {\mathbb R}^N$ be an open bounded set with $C^1$ boundary
and let $Z$ be a vector field of class $C^1(\Omega)\cap
C(\overline{\Omega})$. Then,
\[ \int_\Omega {\rm div}Z(y)\,dy=
\int_{\partial\Omega}\langle Z(x), \nu(x)\rangle \,d\sigma(x).\]
\end{prop}

Here, $\nu(x)$ is the unit outward normal to $\partial\Omega$ at $x$ and $%
d\sigma$ denotes the Euclidian measure on $\partial\Omega$. We denote by $%
\langle \cdot, \cdot \rangle $ the canonical inner product on ${\mathbb R}^N\times
{\mathbb R}^N$.

If $u$ is a differentiable function defined near $\partial\Omega$, we can
define the \emph{normal derivative} of $u$ on $\partial \Omega$ by 
\begin{equation*}
\frac{\partial u}{\partial \nu}=\langle \nabla u, \nu \rangle,\quad %
\mbox{where}\ \nabla u=\mathrm{grad}\,u=\left( \frac{\partial u}{\partial x_1%
}, \frac{\partial u}{\partial x_2},\ldots, \frac{\partial u}{\partial x_N}
\right) .
\end{equation*}

If $\Omega$ is a domain for which the divergence theorem applies, then we
have

\begin{prop}[Green's first identity]
\label{green} Assume that $u,v\in C^2(\Omega)\cap
C^1(\overline{\Omega})$. The following holds
\[ \int_\Omega v(x)\Delta u(x)\,dx+\int_\Omega \langle
\nabla u(x),\nabla v(x) \rangle\,dx=\int_{\partial\Omega} v(x)\frac{\partial u}
{\partial \nu}(x)\,d\sigma(x).  \]
\end{prop}

Let $\|\cdot \|_{L^m(\Omega)}$ denote the usual norm on $L^m(\Omega)$, i.e., 
\begin{equation*}
\|u \|_{L^m(\Omega)}= \left(\int_\Omega |u(x)|^m\,dx \right)^{1/m},\quad %
\mbox{if}\ u\in L^m(\Omega)\ \mbox{with}\ 1\leq m<\infty
\end{equation*}
respectively 
\begin{equation*}
\|u\|_{L^\infty(\Omega)}=\inf\{C>0:\ |u(x)|\leq C \ \mbox{a.e. on}\
\Omega\},\quad \mbox{if}\ u\in L^\infty(\Omega).
\end{equation*}

By $W^{1,m}(\Omega)$, $1\leq m\leq \infty$, we understand the Sobolev space
defined by 
\begin{equation*}
W^{1,m}(\Omega)=\left\{ u\in L^m(\Omega)\left| \begin{aligned} & \exists
g_1,g_2,\ldots g_N\in L^m(\Omega)\ \mbox{such that} \\ & \int_{\Omega}
u\frac{\partial \phi}{\partial x_i}=-\int_{\Omega}g_i \phi, \ \ \forall
\phi\in C^\infty_c(\Omega),\ \ \forall i=\overline{1,N} \\ \end{aligned} %
\right. \right\}.
\end{equation*}

For $u\in W^{1,m}(\Omega)$ we define $g_i=\frac{\partial u}{\partial x_i}$
and we write 
\begin{equation*}
\nabla u=\mathrm{grad}\,u=\left( \frac{\partial u}{\partial x_1}, \frac{%
\partial u}{\partial x_2},\ldots, \frac{\partial u}{\partial x_N} \right).
\end{equation*}
The Sobolev space $W^{1,m}(\Omega)$ is endowed with the norm 
\begin{equation*}
\|u\|_{W^{1,m}(\Omega)}=\|u\|_{L^m(\Omega)}+\sum_{i=1}^N \left\| \frac{%
\partial u}{\partial x_i}\right\|_{L^m(\Omega)}.
\end{equation*}

For $x\in {\mathbb R}^N$ and $r>0$, set $B_r(x)=\{y\in {\mathbb R}^N:\ |x-y|<r\}$, where $%
|x|=\langle x,x \rangle^{1/2}$.

Let $E(x)$ define the fundamental solution of Laplace's equation $\Delta
E(x)=0$ in ${\mathbb R}^N$ ($N\geq 2$), i.e., 
\begin{equation*}
E(x)=\left\{ 
\begin{array}{lll}
& \displaystyle\frac{1}{2\pi}\,\ln |x|, \quad & x\not=0\ (\mbox{if}\ N=2) \\ 
& \displaystyle\frac{1}{(2-N)\omega_N|x|^{N-2}}, \quad & x\not=0\ (\mbox{if}%
\ N\geq 3)
\end{array}
\right.
\end{equation*}
where $\omega_N$ stands for the area of the unit sphere in ${\mathbb R}^N$. By 
\cite[Proposition 0.7]{3}, we know that the value of $\omega_N$ is 
\begin{equation*}
\omega_N=\frac{2\pi^{N/2}} {\Gamma(N/2)}
\end{equation*}
where $\Gamma(s)$ represents the Gamma function defined for $\mathrm{Re}\,
s>0$ by 
\begin{equation*}
\Gamma(s)=\int_0^\infty e^{-t}t^{s-1}\,dt.
\end{equation*}

Let $\Omega\subset {\mathbb R}^N$ be an open, bounded subset with $C^2$ boundary. For
a continuous function $h$ on $\partial\Omega$, the \emph{double layer
potential with moment} $h$ is defined as 
\begin{equation}
\bar{\bar u}_h(y)= \int_{\partial\Omega} h(x)\frac{\partial E}{\partial \nu}
(x-y)\,d\sigma(x).  \label{defi}
\end{equation}

For details about the next results, we refer to \cite{3}.

\begin{prop} \label{es} 
If $h$ is a continuous function on $\partial\Omega$, then
\begin{itemize}
\item[(a)] $\bar{\bar u}_h(y)$ is well defined for all $y\in {\mathbb R}^N$.
\item[(b)] $\Delta \bar{\bar u}_h(y)=0$ for all $y\not\in \partial\Omega$.
\end{itemize}
\end{prop}

\begin{lem}[Gauss' Lemma] \label{Gauss}
Let $\bar{\bar{v}}$ be the double layer potential
with moment $h\equiv 1$, i.e.,
\[ \bar{\bar{v}}(y)=\int_{\partial\Omega}\frac{\partial E}{\partial\nu}(x-y)\,d\sigma(x) .\]
Then, we have
\[ \bar{\bar{v}}(y)=
\left\{\begin{array}{lll}
& 1\quad & \mbox{if}\ y\in \Omega,\\
& 1/2 \quad & \mbox{if}\ y\in \partial\Omega,\\
& 0\quad & \mbox{if}\ y\in {\mathbb R}^N\setminus \overline{\Omega}.
\end{array} \right.\]
\end{lem}

The next result states the limits of the $\bar{\bar{u}}_h(y)$ (defined by (%
\ref{defi})) as we approach $\partial\Omega$ from the interior or exterior
of $\Omega$.

\begin{prop} \label{as} 
Let $h$ be continuous on $\partial \Omega$ and
$y_0\in \partial\Omega$. Then,
\begin{equation}
\lim_{\Omega\ni y\to y_0}\bar{\bar{u}}_h(y)=\frac{1}{2}h(y_0)+
\bar{\bar{u}}_h(y_0)
\ \mbox{and}\
\lim_{{\mathbb R}^N\setminus{\overline{\Omega}}\ni y\to y_0}
\bar{\bar{u}}_h(y)=-\frac{1}{2}h(y_0)+\bar{\bar{u}}_h(y_0).
\label{ie} \end{equation}
\end{prop}

\begin{rem} \label{integ} If $h\in C(\partial\Omega)$ then
$\bar{\bar{u}}_h\in C(\partial\Omega)\cap
L^m(\Omega)$, for each $1\leq m\leq \infty$.
\end{rem}

Indeed, by Propositions~\ref{es} and \ref{as}, the function $%
\phi:\overline\Omega\to {\mathbb R}$ defined by $\phi(y)=\bar{\bar{u}}_h(y)$, $%
\forall y\in \Omega$ and $\phi(y_0)=\frac{1}{2}h(y_0)+\bar{\bar{u}}_h(y_0)$, 
$\forall y_0\in \partial\Omega$ is continuous on $\overline\Omega $. It
follows that ${\bar{\bar u}_h}\in C(\partial\Omega)$ and $\phi\in
L^\infty(\Omega)$. But $\phi\equiv \bar{\bar u}_h$ on $\Omega$ so that $\bar{%
\bar{u}}_h\in L^\infty(\Omega)$. Thus, for each $1\leq m<\infty$, we have 
\begin{equation*}
\int_\Omega |\bar{\bar u}_h|^m\,dx\leq \|\bar{\bar u}_h\|
^m_{L^\infty(\Omega)}\,\mathrm{meas}\,(\Omega)<\infty,
\end{equation*}
which shows that $\bar{\bar u}_h\in L^m(\Omega)$.

\section{Main results}

Let $\Omega\subset {\mathbb R}^N$ be an open bounded set with smooth boundary and $%
A=(a_i)_{i\in I}$ be a finite family of points in $\Omega$.

We assume throughout that $f\in C(\overline{\Omega})\cap C^1(\Omega\setminus
A)$ and, for some $\alpha \in (0,1)$, 
\begin{equation}
\limsup\limits_{x\to a_i}\frac{|f(x)-f(a_i)|}{|x-a_i|^\alpha}<\infty,\quad
\forall i\in I.  \tag{$H$}
\end{equation}

We adopt the following notations 
\begin{equation*}
\oint_{\Omega}f\,dx=\frac{1}{\mathrm{meas}\,(\Omega)}\int_\Omega f(x)\,dx \ %
\mbox{and} \ \oint_{\partial \Omega}f\,d\sigma(x)=\frac{1}{\mathrm{meas}%
\,(\partial \Omega)}\int_{\partial \Omega} f(x)\,d\sigma(x).
\end{equation*}

\begin{th1} \label{teo1}
Suppose $f\in W^{1,p}(\Omega)$ for some $p\in (N,\infty]$. Then
%the exact value of $f(y)$ with respect to
%$\bar{\bar u}_f(y)$ ($y\in \Omega$) is given by
\begin{equation} \label{f1}
f(y)=\bar{\bar u}_f(y)
-\int_{\Omega} \langle \nabla E(x-y),\nabla f(x)\rangle\,dx,\qquad
\forall y\in \Omega
\end{equation}
resp.,
\begin{equation} \label{fig}
\int_\Omega f(x)\,dx=\frac{1}{N}\int_{\partial\Omega} f(x)\langle
x-y,\nu \rangle\,d\sigma(x)-\frac{1}{N}\int_{\Omega}\langle \nabla f(x),
x-y \rangle \,dx,\ \ \forall y\in {\mathbb R}^N.
\end{equation}
\end{th1}

\begin{proof}%[Proof of Theorem~\ref{teo1}] 
Let $y\in \Omega$ be fixed. 
We first recall that, for each $\gamma\in (0,N)$, the mapping
$x\longmapsto |x-y|^{-\gamma}\in L^1(\Omega)$. 
Indeed, for $r>0$ fixed so that $B_r(y)\subset\subset \Omega$,
we have
\[\begin{split}
\int_\Omega \frac{dx}{|x-y|^\gamma}&=\int_{\Omega\setminus
B_r(y)}\frac{dx}{|x-y|^\gamma}+\int_{B_r(y)}\frac{dx}{|x-y|^\gamma}\\
&\leq \frac{{\rm meas}\,(\Omega)}{r^\gamma}+
\lim_{\epsilon\to 0}\int_\epsilon^r\left(\int_{\partial B_\rho(y)}\frac{d\sigma(x)}
{|x-y|^\gamma}\right)d\rho\\
&= \frac{{\rm meas}\,(\Omega)}{r^\gamma}+
\lim_{\epsilon\to 0}\int_\epsilon^r \frac{{\rm meas}\,(\partial B_\rho(y))}{\rho^\gamma}
\,d\rho\\
&= \frac{{\rm meas}\,(\Omega)}{r^\gamma}+\frac{\omega_N r^{N-\gamma}}{N-\gamma}<\infty.
\end{split} \]
We now define $F:\overline{\Omega}\setminus \{y\}\to {\mathbb R}^N$
as follows
\[ F(x)=(f(x)-f(y))\nabla E(x-y)=\frac{f(x)-f(y)}{\omega_N|x-y|^N}\,(x-y).\]
Note that $F(x)$ is not smooth for all $x\in \Omega$. We overcome this problem
by choosing $\epsilon>0$ small enough such that
$B_\epsilon(y)$ resp., $B_\epsilon(a_i)$ ($a_i\in A\setminus\{y\}$) is contained
within $\Omega$ and each two such balls are disjoint. Therefore,
$F\in C^1(D_\epsilon)\cap C(\overline{D}_\epsilon)$ where $D_\epsilon=\Omega\setminus
\left(\cup_{i\in I}\overline{B}_\epsilon(a_i)
\cup\overline{B}_\epsilon(y)\right)$. Using the Divergence Theorem, we arrive at
\begin{equation}
\begin{split}
\int_{D_\epsilon}{\rm div}\,F(x)\,dx&=
\int_{\partial\Omega}(f(x)-f(y))\,\frac{\partial E}{\partial \nu}
(x-y)\,d\sigma(x)\\
&\quad -\frac{1}{\omega_N\epsilon^{N-1-\alpha}}\int_{\partial B_\epsilon(y)}
\frac{f(x)-f(y)}{|x-y|^\alpha}\,d\sigma(x) \\
& \quad -\frac{1}{\omega_N}\sum_{i\in I, a_i\not=y}
\int_{\partial B_\epsilon(a_i)}
\frac{f(x)-f(y)}{\epsilon|x-y|^N}\langle x-y,x-a_i\rangle\, d\sigma(x). \label{ct0}
\end{split} \end{equation}
We see that
\begin{equation} \lim_{\epsilon\to 0}\frac{1}{\epsilon^{N-1-\alpha}}
\int_{\partial B_\epsilon(y)}\frac{f(x)-f(y)}{|x-y|^\alpha}\,d\sigma(x)=0.
\label{ct1} \end{equation}
Indeed, in view of $(H)$,
for some constant $L>0$ and $\epsilon>0$ small enough, we have
\[ \begin{split}
0 &\leq \frac{1}{\epsilon^{N-1-\alpha}}
\left|\int_{\partial B_\epsilon(y)}
\frac{f(x)-f(y)}{|x-y|^\alpha}\,
d\sigma(x)\right|\\
& \leq \frac{L}{\epsilon^{N-1-\alpha}}
\int_{\partial B_\epsilon(y)}d\sigma(x)=L\omega_N \epsilon^\alpha\to 0\ \mbox{as}\ \epsilon\to 0.
\end{split}\]
Notice that, for each
$i\in I$ with $a_i\not=y$, there exists
a constant $C_i>0$ such that
\[ |f(x)-f(y)|\leq C_i|x-y|^{N-1},\ \forall x\in \overline{B}_\epsilon(a_i)\]
(since $y\not \in \overline{B}_\epsilon(a_i)$). Hence
\begin{equation}
\begin{split}
\left|\int_{\partial B_\epsilon(a_i)}
\frac{f(x)-f(y)}{\epsilon|x-y|^N}\langle x-y,x-a_i\rangle\, d\sigma(x)
\right| & \leq \int_{\partial B_\epsilon(a_i)}\frac{|f(x)-f(y)|}{|x-y|^{N-1}}\,
d\sigma(x)\\
& \leq C_i\omega_N \epsilon^{N-1}\to 0\quad \mbox{as}\ \epsilon\to 0,
\label{ct2}
\end{split}
\end{equation}
provided $i\in I$ such that $a_i\not=y$.
By (\ref{ct0})--(\ref{ct2}), it follows that
\begin{equation}
\begin{split}
\lim_{\epsilon\to 0}\int_{D_\epsilon}{\rm div}\,F(x)\,dx &=
\int_{\partial\Omega}
(f(x)-f(y))\frac{\partial E}{\partial \nu}(x-y)\,d\sigma(x)\\
&=\int_{\partial\Omega}
f(x)\frac{\partial E}{\partial \nu}(x-y)\,d\sigma(x)-f(y) \label{aux}
\end{split} \end{equation}
by using Gauss' Lemma.
On the other hand, for each $x\in D_\epsilon$,
\[ \begin{split}
{\rm div}\, F(x)&=\langle \nabla f(x), \nabla E(x-y)\rangle+
(f(x)-f(y))\Delta_x E(x-y)\\
& = \langle \nabla f(x), \nabla E(x-y)\rangle
\end{split} \]
since $x\longmapsto E(x-y)$ is harmonic on ${\mathbb R}^N\setminus \{y\}$.
By H\"older's inequality, we obtain
\[ \int_{\Omega}|\langle \nabla f(x), \nabla E(x-y)\rangle|\,dx\leq
\frac{\|\nabla f\|_{L^p(\Omega)}}{\omega_N}\left(\int_{\Omega}\frac{dx}
{|x-y|^{(N-1)p'}}\right)^{\frac{1}{p'}}<\infty \]
which is due to $|\nabla f|\in L^p(\Omega)$ and $(N-1)p'<N$.
Hence, the mapping $x\longmapsto \langle \nabla f(x),\nabla E(x-y)
\rangle $ is integrable on $\Omega$. Thus, using (\ref{aux}) 
we deduce that
\[ \begin{split}
\int_{\Omega}\langle \nabla f(x), \nabla E(x-y)\rangle \,dx &=
\lim_{\epsilon\to 0}
\int_{D_\epsilon}\langle \nabla f(x), \nabla E(x-y)\rangle \,dx\\
&=\int_{\partial\Omega}
f(x)\frac{\partial E}{\partial \nu}(x-y)\,d\sigma(x)-f(y)
\end{split} \]
which concludes our first assertion.

Let $y\in {\mathbb R}^N$ be arbitrary.
We define $G:\overline{\Omega}\to {\mathbb R}^N$ by $G(x)=f(x)(x-y)$.
Let $\epsilon>0$ be small such that
$\overline{B}_\epsilon(a_i)\subset \Omega$,
$\forall i\in I$ and $\overline{B}_\epsilon(a_i)\cap
\overline{B}_\epsilon(a_j)=\emptyset$, $\forall i,j\in I$
with $i\not=j$.
Set $U_\epsilon=\Omega\setminus \left(\cup_{i\in I}
\overline{B}_\epsilon(a_i)\right)$. We have
$G\in C^1(U_\epsilon)\cap C(\overline{U}_\epsilon)$. By Proposition~\ref{div},
we find that 
\begin{equation} \label{char1}
\begin{split}
\int_{U_\epsilon}{\rm div}\, G(x)\,dx &=\int_{\partial\Omega}f(x)\langle
x-y,\nu \rangle \,d\sigma(x)\\
& -\sum_{i\in I}\int_{\partial B_\epsilon(a_i)}
\frac{f(x)}{\epsilon}\langle x-y, x-a_i\rangle\,d\sigma(x).\\
\end{split}
\end{equation}
For each $i\in I$, we have
\[\begin{aligned}
\left|
\int_{\partial B_\epsilon(a_i)}\frac{f(x)}{\epsilon}
\langle x-y,x-a_i \rangle \,d\sigma(x)
\right|& \leq
\int_{\partial B_\epsilon(a_i)}
\frac{|f(x)|}{\epsilon}|\langle x-y,x-a_i \rangle| \,d\sigma(x)\\
& \leq
\int_{\partial B_\epsilon(a_i)}
|f(x)||x-y|\,d\sigma(x)\\
& \leq C_i \|f\|_{L^\infty(\Omega)}{\rm meas}\,(\partial B_\epsilon(a_i))\\
& =C_i \|f\|_{L^\infty(\Omega)}\omega_N \epsilon^{N-1}\to 0\ \mbox{as}\
\epsilon \to 0\\
\end{aligned} \]
for some constant $C_i>0$ that satisfies
$|x-y|\leq C_i$, $\forall x\in \partial B_k(a_i)$, $\forall k\in (0,\epsilon]$.

It follows that
\begin{equation}\label{char2}
\lim_{n\to \infty}
\int_{\partial B_\epsilon(a_i)}\frac{f(x)}{\epsilon}
\langle x-y,x-a_i \rangle \,d\sigma(x)=0, \quad \forall i\in I.
\end{equation}
We see that
\[ {\rm div}\, G(x)=\langle \nabla f(x),x-y\rangle+N f(x),\quad
\forall x\in U_\epsilon. \]
By $f\in C(\overline{\Omega})\cap W^{1,p}(\Omega)$ and H\"{o}lder's
inequality, we deduce 
$f\in L^1(\Omega)$ and
\[ \begin{aligned}
\int_{\Omega} |\langle \nabla f(x), x-y \rangle |\,dx & \leq
\int_\Omega |\nabla f(x)||x-y|\,dx\\
& \leq \left(\int_{\Omega} |\nabla f(x)|^p\,dx\right)^{\frac{1}{p}}
\left(\int_\Omega |x-y|^{p'}\,dx\right)^{\frac{1}{p'}}\\
&=\|\nabla f\|_{L^p(\Omega)}
\left(\int_\Omega |x-y|^{p'}\,dx\right)^{\frac{1}{p'}}<\infty.
\end{aligned}
\]
Therefore,
\begin{equation} \label{char3}
\lim_{\epsilon\to 0}\int_{U_\epsilon} {\rm div}\,G(x)\,dx=
\int_{\Omega} \langle \nabla f(x),x-y\rangle +N \int_{\Omega} f(x)\,dx.
\end{equation}
Passing to the limit $\epsilon\to 0$ in (\ref{char1}) and using (\ref{char2})
resp., (\ref{char3}), we conclude that
\[
\int_{\Omega} \langle \nabla f(x),x-y\rangle+
N \int_{\Omega} f(x)\,dx=
\int_{\partial \Omega} f(x)\langle x-y,\nu \,d\sigma(x)\]
which proves (\ref{fig}).
\end{proof}

To our next aim, we recall the following results.

\begin{lem} \label{conv} Let $\Omega\subset {\mathbb R}^N$
be an open set. Let $(h_n)$ be a sequence in $ L^p(\Omega)$, $1\leq p\leq \infty$, and let
$h\in L^p(\Omega)$ be such that $\|h_n-h\|_{L^p(\Omega)}\to 0$.

Then, there exists a subsequence $(h_{n_k})$ and a function $\phi\in L^p(\Omega)$
such that

\begin{itemize}
\item[(a)] $h_{n_k}(x)\to h(x)$ a.e. in $\Omega$,
\item[(b)] $|h_{n_k}(x)|\leq \phi(x)$ $\forall k$, a.e. in $\Omega$.
\end{itemize}
\end{lem}
The interested reader may find the proof of Lemma~\ref{conv} in 
\cite[Theorem IV.9]{1}.

\begin{lem} \label{dens} Suppose that $\Omega$ is of class $C^1$ and
let $u\in W^{1,p}(\Omega)$ with $1\leq p<\infty$.

Then, there exists a sequence
$(u_n)$ in $C^\infty_c({\mathbb R}^N)$ such that $u_n|_{\Omega}\to u$
in $W^{1,p}(\Omega)$. In other words, the restrictions to $\Omega$
of functions belonging to $C^\infty_c({\mathbb R}^N)$ form a subspace
which is dense in $W^{1,p}(\Omega)$.
\end{lem}
For the proof of Lemma~\ref{dens} we refer to \cite[Corollary IX.8]{1}.

We are now ready to give a representation theorem of functions in any
Sobolev space $W^{1,p}(\Omega)$, $p\in (N,\infty)$. More precisely, we prove

\begin{th1} \label{teo2}
Let $\Omega$ be an open bounded $C^1$ set in ${\mathbb R}^N$, $N\geq 2$.
Then, for any $g\in W^{1,p}(\Omega)$ with
$p\in (N,\infty)$, there exists a sequence $(g_n)\subset C^\infty_c
({\mathbb R}^N)$ so that
\begin{equation}
\begin{split}
g(y)& =\lim_{n\to \infty}\int_{\partial\Omega}g_n(x)\frac{\partial E}
{\partial \nu}(x-y)\,d\sigma(x)\\
& \quad -\int_\Omega \langle
\nabla E(x-y),\nabla g(x)\rangle \,dx\quad \mbox{a.e. } y\in \Omega.
\end{split} \end{equation}
\end{th1}

\begin{proof}
By Lemma~\ref{dens}, we know that there exists a sequence
$g_n\in C^\infty_c({\mathbb R}^N)$ such that $g_n|_{\Omega}\to g$ in
$W^{1,p}(\Omega)$. Hence,
\[\lim_{n\to \infty} \|g_n|_{\Omega}-g\|_{L^p(\Omega)}=0\ \
\mbox{and}\ \ \lim_{n\to \infty} \left\| \frac{\partial g_n}
{\partial x_i}-\frac{\partial g}{\partial x_i}\right\|_{L^p(\Omega)}=0,\
\forall i=\overline{1,N}.\]
Applying Lemma~\ref{conv} we have that, up to a subsequence
(relabelled $(g_n)$),
\begin{equation} \label{point}
g_n|_{\Omega}\to g\quad \mbox{a.e. in}\ \Omega.
\end{equation}
Using Theorem~\ref{teo1}, we obtain
\begin{equation} \label{st}
g_n(y)=\int_{\partial\Omega}g_n(x)\frac{\partial E}{\partial \nu}(x-y)\,
d\sigma(x)-\int_{\Omega}\langle \nabla E(x-y), \nabla g_n(x)\rangle \,dx,
\quad \forall y\in \Omega.\end{equation}
We now show that
\begin{equation} \label{st2}
\lim_{n\to \infty}
\int_{\Omega}\langle \nabla E(x-y), \nabla g_n(x)\rangle \,dx=
\int_{\Omega} \langle \nabla E(x-y), \nabla g(x)\rangle \,dx,\quad
\forall y\in \Omega.
\end{equation}
Indeed, by H\" {o}lder's inequality, we deduce
$$ \begin{aligned}
0 & \leq \int_{\Omega} |\langle E(x-y),\nabla g_n(x)-\nabla g(x)
\rangle | \,dx\\
& =\int_{\Omega} \left| \sum_{i=1}^N
\frac{\partial E}{\partial x_i}(x-y)
\,\frac{\partial(g_n-g)}{\partial x_i} \right| 
\leq \sum_{i=1}^N \int_{\Omega} \left|
\frac{\partial E}{\partial x_i}(x-y)
\,\frac{\partial(g_n-g)}{\partial x_i}
\right| \,dx\\
& \leq \sum_{i=1}^N
\left(\int_{\Omega}\left|
\frac{\partial E}{\partial x_i}(x-y)\right|^{p'}\,dx\right)^{1/p'}
\cdot \left(\int_{\Omega}\left|
\frac{\partial(g_n-g)}{\partial x_i}\right|^p\,dx\right)^{1/p}\\
& \leq
\left(\int_{\Omega}\left|
\nabla E(x-y)\right|^{p'}\,dx\right)^{1/p'}
\sum_{i=1}^N \left\|\frac{\partial(g_n-g)}{\partial x_i}
\right\|_{L^p(\Omega)}
\\
& \leq \frac{1}{\omega_N}\left(\int_{\Omega}\frac{dx}{|x-y|^{(N-1)p'}}
\right)^{1/p'}\cdot
\sum_{i=1}^N \left\|\frac{\partial(g_n-g)}{\partial x_i}
\right\|_{L^p(\Omega)}\to 0 \ \mbox{as}\ n\to \infty\\
\end{aligned}
$$
By (\ref{point})--(\ref{st2}) we conclude the proof.
\end{proof}

\section{Special cases}

A function $u\in C^2(\Omega)$ is called \emph{harmonic} in $\Omega$ if it
satisfies $\Delta u=0$ in $\Omega$.

The mean value theorem for harmonic functions says that the function value
at the center of the ball $B_R(a)\subset \Omega$ is equal to the integral
mean values over both the surface $\partial B_R(a)$ and $B_R(a)$ itself.
More precisely,

\begin{prop}[Theorem~2.1 in \cite{4}] \label{trud}
Let $u\in C^2(\Omega)\cap C(\overline{\Omega})$ satisfy $\Delta u=0$ in
$\Omega$. Then for any ball $B_R(a)\subset \Omega$, we have
\begin{equation} \label{q1}
u(a)=\oint_{\partial B_R(a)}u(x)\,d\sigma(x),
\end{equation}
\begin{equation} \label{q2}
u(a)=\oint_{B_R(a)}u(x)\,dx.
\end{equation}
\end{prop}

The \emph{Poisson integral formula}, together with an approximation
argument, gives the representation form for harmonic functions $u\in
C^2(B_R(a))\cap C(\overline{B}_R(a))$, that is (see \cite{4}, pp. 20) 
\begin{equation}  \label{poisson}
u(y)= \frac{R^2-|y-a|^2}{R\omega_N}\int_{\partial B_R(a)}\frac{u(x)}{|x-y|^N}%
\, d\sigma(x),\quad \forall y\in B_R(a).
\end{equation}
Moreover, we have

\begin{prop}[Theorem 2.6 in \cite{4}] \label{diric}
Let $\varphi$ be a continuous function on $\partial B$. Then the function
$u$ defined by
\begin{equation} \label{pi}
u(y)=\left\{ \begin{aligned}
& \frac{R^2-|y-a|^2}{R\omega_N}\int_{\partial B_R(a)}\frac{u(x)}{|x-y|^N}\,
d\sigma(x),\quad \forall y\in B_R(a),\\
& \varphi(y), \quad \forall y\in \partial B_R(a)
\end{aligned} \right.
\end{equation}
belongs to $C^2(B_R(a))\cap C(\overline{B}_R(a))$ and satisfies
$\Delta u=0$ in $B_R(a)$.
\end{prop}

It is now natural to ask what are the corresponding representation formulas
for functions satisfying weaker regularity assumptions and not necessarily
harmonic.

To this aim, we state some consequences of Theorem~\ref{teo1}, whose
preliminary assumptions are self-understood. As a common hypothesis for
Corollaries~\ref{mean}--\ref{second}, we have $f\in W^{1,p}(\Omega)$ for
some $p\in (N,\infty]$.

\begin{cor} \label{mean}
For any ball $B_R(a)\subset \Omega$, we have
\begin{equation} \label{mat}
f(y)=
\int_{\partial B_R(a)}\frac{f(x)\langle
x-y, x-a \rangle}{R\omega_N|x-y|^N}\,d\sigma(x)
-\int_{B_R(a)}\frac{\langle \nabla f(x),x-y\rangle}{\omega_N|x-y|^N}\,dx,
\end{equation}
where $y\in B_R(a)$ is arbitrary.
\end{cor}

Using Proposition~\ref{diric} and Corollary~\ref{mean}, we arrive at

\begin{cor} \label{oi} For any $a\in \Omega$ and $R>0$ such that
$B_R(a)\subset \Omega$, we find
\begin{equation} \label{com}
\begin{aligned} f(y)=& \chi(y)+
\int_{\partial B_R(a)}\frac{\langle
y-a, y-x \rangle}{R\omega_N|x-y|^N}\,f(x)\,d\sigma(x)\\
& -\int_{B_R(a)}\frac{\langle \nabla f(x),x-y\rangle}{\omega_N|x-y|^N}\,dx,
\quad \forall y\in B_R(a)
\end{aligned}
\end{equation}
where $\chi$ is the unique classical solution of the Dirichlet problem
\[ \left\{ \begin{aligned}
& \Delta u=0, \quad \mbox{in}\ B_R(a)\\
& u=f, \quad \mbox{on}\ \partial B_R(a).
\end{aligned} \right.
\]
\end{cor}

\begin{cor} \label{rep}
The following representation formula holds
\begin{equation} \label{rp0}
\begin{aligned}
f(y)=&\oint_{\Omega} f(x)\,dx+\int_{\partial\Omega}
\left(\frac{\langle x-y,\nu\rangle}{\omega_N |x-y|^N}-\frac{
\langle x-z,\nu\rangle}{N\,{\rm meas}\,(\Omega)}
\right)f(x)\,d\sigma(x)\\
&-\int_{\Omega}
\left(
\frac{\langle \nabla f(x),x-y\rangle}
{\omega_N |x-y|^N}-\frac{\langle
\nabla f(x),x-z\rangle}{N{\rm meas}\,(\Omega)}\right)
\,dx, \quad \forall y\in \Omega,\ \ \forall z\in {\mathbb R}^N.\\
\end{aligned}
\end{equation}
In particular, for $z=y$ we obtain
\begin{equation} \label{rp1}
\begin{aligned}
f(y)=&\oint_{\Omega} f(x)\,dx+\int_{\partial\Omega}
\left(\frac{1}{\omega_N |x-y|^N}-\frac{1}{N{\rm meas}\,(\Omega)}\right)
f(x)\langle x-y,\nu\rangle \,d\sigma(x)\\
&-\int_{\Omega}
\left(\frac{1}{\omega_N |x-y|^N}-\frac{1}{N{\rm meas}\,(\Omega)}\right)
\langle \nabla f(x),x-y\rangle \,dx, \quad \forall y\in \Omega.\\
\end{aligned}
\end{equation}
\end{cor}

\begin{cor} \label{cerc}
For each $a\in \Omega$ and $R>0$ such that $B_R(a)\subset
\Omega$, we obtain
%\begin{equation} \label{repr}
\[\begin{aligned}
f(y)=&\oint_{B_R(a)}f(x)\,dx-\oint_{\partial B_R(a)}f(x)\,d\sigma(x)+
\int_{\partial B_R(a)} \frac{f(x)\langle x-y,x-a \rangle}
{R\omega_N|x-y|^N}\,d\sigma(x)\\
&-\frac{1}{\omega_N}\int_{B_R(a)}\left(\frac{\langle \nabla f(x),x-y
\rangle}{|x-y|^N}-\frac{\langle \nabla f(x),x-a \rangle}{R^N}\right)\,dx,\qquad
\forall y\in B_R(a).\end{aligned} \]
%\end{equation}
The particular case $y=a$ leads to
\begin{equation} \label{rep3}
f(a)=\oint_{B_R(a)}f(x)\,dx -\frac{1}{\omega_N}\int_{B_R(a)}
\left(\frac{1}{|x-a|^N}-\frac{1}{R^N}\right)\langle \nabla f(x),x-a\rangle\,dx.
\end{equation}
resp.,
\begin{equation} \label{rep2}
f(a)=\oint_{\partial B_R(a)}f(x)\,d\sigma(x)-\frac{1}{\omega_N}
\int_{B_R(a)}\frac{\langle\nabla f(x),x-a \rangle }{|x-a|^N}\,dx.
\end{equation}
\end{cor}

\begin{cor} \label{c1}
An arbitrary value of $f$ is below compared with
the double layer potential
with moment $f$
\begin{equation} \label{sec} 
\left|f(y)-\bar{\bar u}_f(y)\right|\leq
\frac{\|\nabla f\|_{L^p(\Omega)}}
{\omega_N}\left(\int_\Omega \frac{dx}{|x-y|^{(N-1)p'}}\right)^{\frac{1}{p'}},
\quad \forall y\in \Omega
\end{equation}
where $p'$ denotes the conjugate coefficient of $p$ (i.e., $1/p+
1/p'=1$).
Moreover, for $y\in \Omega$ fixed,
the equality in (\ref{sec}) is established for the nontrivial function
$f(x)=\pm |x-y|$ if $p=\infty$ resp., $f(x)=\pm |x-y|^\beta$ with
$\beta=(p-N)/(p-1)$ if $p\in (N,\infty)$.
\end{cor}

\begin{proof}
By (\ref{f1}) and H\" {o}lder's inequality, we have
\[ \begin{aligned}
|f(y)-\bar{\bar u}_f(y)| &=\left|\int_{\Omega} \langle
\nabla E(x-y),\nabla f(x)\rangle\,dx\right|=\left|
\int_\Omega \frac{\langle x-y, \nabla f(x)\rangle }{\omega_N |x-y|^N}\,dx
\right|\\
& \leq \frac{1}{\omega_N}\int_{\Omega}\frac{|\langle
x-y, \nabla f(x)\rangle|}{|x-y|^N}\,dx\leq \frac{1}{\omega_N}
\int_{\Omega}\frac{|\nabla f(x)|}{|x-y|^{N-1}}\,dx\\
& \leq \frac{1}{\omega_N} \left(\int_{\Omega}|\nabla f(x)|^p\,dx
\right)^{1/p}
\left(\int_{\Omega}\frac{dx}{|x-y|^{(N-1)p'}}\right)^{1/p'}\\
&=\frac{\|\nabla f\|_{L^p(\Omega)}}{\omega_N}
\left(\int_{\Omega}\frac{dx}{|x-y|^{(N-1)p'}}\right)^{1/p'}.\\
\end{aligned} \]
Let $y\in \Omega$ be fixed. We define $f^{\pm}_{p,y}:\overline{\Omega}\to
{\mathbb R}$ by
\[ f_{p,y}^\pm(x)=\left\{
\begin{aligned}
& \pm |x-y|,\quad \mbox{if}\ p=\infty\\
& \pm |x-y|^{\frac{p-N}{p-1}}, \quad \mbox{if}\ p\in (N,\infty).
\end{aligned} \right.\]
Clearly, we have $f_{p,y}^\pm \in C(\overline{\Omega})$. Moreover,
$f_{p,y}^\pm \in C^1(\Omega\setminus \{y\})$ and
\begin{equation} \label{pw}
\nabla f_{p,y}^\pm(x)=\left\{
\begin{aligned}
& \pm \frac{x-y}{|x-y|},\quad \forall x\in \Omega\setminus \{y\},
\quad \mbox{if}\ p=\infty\\
& \pm \frac{p-N}{p-1}\frac{x-y}{|x-y|^{\frac{p+N-2}{p-1}}},\quad
\forall x\in \Omega\setminus \{y\},
\quad \mbox{if}\ p\in (N,\infty).
\end{aligned} \right. \end{equation}
Since $C(\overline \Omega)\subset L^p(\Omega)$, we infer that
$f_{p,y}^\pm \in W^{1,p}(\Omega)$ and
\[ \left\|\nabla f_{p,y}^\pm(x)\right\|_{L^p(\Omega)}=\left\{
\begin{aligned}
& 1, \quad \mbox{if}\ p=\infty\\
& \frac{p-N}{p-1}\left(\int_{\Omega}\frac{dx}{|x-y|^{(N-1)p'}}\right)^{1/p},
\quad \mbox{if}\ p\in (N,\infty).
\end{aligned} \right.\]
It follows that the right hand side (RHS) of (\ref{sec}) for $f_{p,y}^\pm$ is
\begin{equation} \label{com1}
{\rm RHS}=
\left\{
\begin{aligned}
& \frac{1}{\omega_N}\left(\int_{\Omega}\frac{dx}{|x-y|^{N-1}}\right),
\quad \mbox{if}\ p=\infty\\
& \frac{p-N}{\omega_N(p-1)}
\int_{\Omega}\frac{dx}{|x-y|^{(N-1)p'}},
\quad \mbox{if}\ p\in (N,\infty).
\end{aligned} 
\right.
\end{equation}
By (\ref{f1}) and (\ref{pw}), we have that the left hand side (LHS) of (\ref{sec})
for $f_{p,y}^\pm $ is
\begin{equation} \label{com2}
\begin{aligned}
{\rm LHS}&=\left|\int_{\Omega}\langle \nabla E(x-y),
\nabla f^\pm_{p,y}(x)\rangle\,dx\right|=\left|
\int_\Omega \frac{\langle x-y, \nabla f_{p,y}^\pm (x)\rangle}{\omega_N |x-y|^N}\,dx
\right|\\
& =
\left\{
\begin{aligned}
& \frac{1}{\omega_N}\left(\int_{\Omega}\frac{dx}{|x-y|^{N-1}}\right),
\quad \mbox{if}\ p=\infty\\
& \frac{p-N}{\omega_N(p-1)}
\int_{\Omega}\frac{dx}{|x-y|^{(N-1)p'}},
\quad \mbox{if}\ p\in (N,\infty).
\end{aligned} 
\right.
\end{aligned}
\end{equation}
Using (\ref{com1}) and (\ref{com2}) we obtain equality in (\ref{sec}) for
$f(x)=f^\pm_{p,y}(x)$.
\end{proof}

\begin{cor} \label{corl2} For $a\in \Omega$ and $R>0$ such that
$B=B_R(a)\subset \overline{B}_R(a)\subset \Omega$, we have
\begin{equation}
\left|f(a)-\oint_{\partial B}f(x)\,d\sigma(x)
\right| \leq 
\omega_N^{\frac{1}{p'}-1}\left(\frac{R^{N-(N-1)p'}}{N-(N-1)p'}\right)^{
\frac{1}{p'}}
\|\nabla f\|_{L^p(B)}\,. \label{eq6}
\end{equation}
%\left|f(y)-\oint_{B}f(x)\,dx
%\right| &\leq &
%\frac{Np'}{N+p'}\,\omega_N^{\frac{1}{p'}-1}\left(\frac{R^{N-(N-1)p'}}{N-(N-1)p'}\right)^{
%\frac{1}{p'}}
%\|\nabla f\|_{L^p(B)}\,. \label{eq7}
%\end{eqnarray}
Moreover, the constant is sharp and the 
function
$f(x)=\pm |x-a|$ if $p=\infty$ resp.,
$f(x)=\pm |x-a|^{(p-N)/(p-1)}$ if $p\in (N,\infty)$
achieves the equality.
\end{cor}

\begin{proof}
Note that $f\in C(\overline{B})\cap C^1(B\setminus A_i)$ resp.,
$f\in W^{1,p}(B)$ with $p\in (N,\infty]$. Therefore, we can apply
Corollary~\ref{c1} with $y=a$ and $\Omega=B$. More precisely,
\begin{equation} \label{tec}
\left|f(a)-\int_{\partial B}f(x)\frac{\partial E}{\partial \nu}(x-a)\,d\sigma(x)
\right|
\leq \frac{\|\nabla f\|_{L^p(B)}}{\omega_N}\left(
\int_B\frac{dx}{|x-a|^{(N-1)p'}}
\right)^{1/p'}
\end{equation}
where the equality holds for
$f(x)=\pm |x-y|$ if $p=\infty$ and $f(x)=\pm |x-y|^{(p-N)/(p-1)}$ 
if $p\in (N,\infty)$.

Notice that, for each $x\in \partial B$, we have
\[ \begin{aligned}
\frac{\partial E}{\partial \nu}(x-a)&=\langle \nabla E(x-a), \nu(x)
\rangle=\langle \frac{x-a}{\omega_N |x-a|^N},\frac{x-a}{|x-a|}
\rangle\\
&=\frac{1}{\omega_N |x-a|^{N-1}}=\frac{1}{\omega_N R^{N-1}}={\rm meas}\, (\partial B).\\ 
\end{aligned}
\]
It follows that 
\begin{equation} \label{nw1}
\int_{\partial B}f(x)\frac{\partial E}{\partial \nu}\,(x-a)\,d\sigma(x)=
\frac{1}{{\rm meas}\,(\partial B)}\int_{\partial B} f(x)\,d\sigma(x)=
\oint_{\partial B} f(x)\,d\sigma(x).
\end{equation}
On the other hand,
\begin{equation} \label{nw2}
\begin{aligned}
\int_B\frac{dx}{|x-a|^{(N-1)p'}}&=\int_0^R \left(
\int_{\partial B_\rho(a)}\frac{d\sigma(x)}{|x-a|^{(N-1)p'}}
\right)\,d\rho\\
&=\int_0^R \left(\frac{1}{\rho^{(N-1)p'}}\int_{\partial B_\rho(a)}
d\sigma(x)\right)\,d\rho
=\int_0^R \frac{\omega_N \rho^{N-1}}{\rho^{(N-1)p'}}\,d\rho\\
&=\frac{\omega_NR^{N-(N-1)p'}}{N-(N-1)p'}.\\
\end{aligned}
\end{equation}
Replacing (\ref{nw1}) and (\ref{nw2}) in (\ref{tec}) we obtain (\ref{eq6}).
\end{proof}

\begin{cor} \label{second}
The following identities hold
\begin{equation}
\begin{aligned}
\int_\Omega \bar{\bar u}_f(y)\,dy 
=& \frac{1}{N}\int_{\partial\Omega} f(x)\langle x-z,\nu \rangle \,d\sigma(x)-
\frac{1}{N}\int_\Omega \langle \nabla f(x),x-z\rangle\,dx\\
& +
\int_\Omega \left(
\int_{\Omega}\langle \nabla E(x-y),\nabla f(x)\rangle\,dx 
\right)dy
,\quad \forall z\in {\mathbb R}^N
\end{aligned}\label{f2}
\end{equation}
resp., 
\begin{equation}
\int_{\partial\Omega}\bar{\bar u}_f(z)\,d\sigma(z)=
\frac{1}{2}\int_{\partial\Omega} f(z)\,d\sigma(z)+
\int_{\partial\Omega} \zeta(z)\,d\sigma(z),
\label{f3}\end{equation}
where we define
\[ \zeta(z)=\lim_{\Omega\ni t\to z}
\int_{\Omega} \langle \nabla E(x-t),\nabla f(x)\rangle \,dx,
\quad \mbox{for each}\ z\in \partial\Omega .\]
\end{cor}

\begin{rem} Note that $\zeta$ is well defined because of (\ref{ie}) and
(\ref{f1}). \end{rem}

\begin{proof}
By virtue of Remark~\ref{integ}, $\bar{\bar{u}}_f\in L^1(\Omega)$.
Obviously, $f\in L^1(\Omega)$ since $f\in C(\overline\Omega)$ and $\Omega$
is bounded. Therefore, we can integrate (\ref{f1}) over $\Omega$ to obtain
\[ \int_\Omega f(y)\,dy=\int_\Omega \bar{\bar u}_f(y)\,dy-
\int_\Omega \left(\int_{\Omega}\langle \nabla E(x-y),\nabla f(x)\rangle\,dx 
\right)dy.\]
Using now (\ref{fig}), we arrive at (\ref{f2}).

Let $z\in \partial\Omega$ be arbitrary. By the continuity of
$f$ on $\overline{\Omega}$ and Proposition 3, we find
\[ \lim_{\Omega\ni y\to z}[f(y)-\bar{\bar{u}}_f(y)]=\frac{f(z)}{2}-
\bar{\bar{u}}_f(z).\]
Combining this with (\ref{f1}), we derive that
\begin{equation} f(z)=2 \bar{\bar{u}}_f(z)-2\zeta(z), \quad \forall z\in \partial\Omega.
\label{sur}\end{equation}
By Remark 1, $\bar{\bar{u}}_f(z)\in C(\partial\Omega)$. Hence 
integrating (\ref{sur}) over $\partial\Omega$ we find (\ref{f3}).
%\[ \oint_{\partial\Omega}f(z)\,d\sigma(z)=2\oint_{\partial\Omega}
%\bar{\bar{u}}_f(z)\,d\sigma(z)-2\oint_{\partial\Omega}\zeta(z)\,d\sigma(z).
%\] This, together with (\ref{f1}), leads to (\ref{f3}).
\end{proof}

\begin{cor}[Gauss' Lemma extension] \label{c2}
Assume $f\in W^{1,p}(\Omega)$, for some $p\in [1,\infty]$. Then
the following representation holds
\begin{equation}
\bar{\bar{u}}_f(y)=\left\{ \begin{aligned}
& f(y)+\int_{\Omega} \langle \nabla E(x-y), \nabla f(x)\rangle\,dx,
\quad \forall y\in \Omega, \ \mbox{if}\ p\in (N,\infty],\\
& \zeta(y)+f(y)/2, \ \ \forall y\in \partial\Omega,\ \mbox{if}
\  p\in (N,\infty],\\
& \int_{\Omega} \langle \nabla E(x-y), \nabla f(x)\rangle\,dx, \quad
\forall y\in {\mathbb R}^N\setminus\overline{\Omega},\ \ \forall p\in [1,\infty]. 
\end{aligned} \right.
\label{potential}
\end{equation}
\end{cor}

\begin{proof} In view of (\ref{f1}) and (\ref{sur}), we need only to show that
\begin{equation} \label{on}
\bar{\bar{u}}_f(y)=
\int_{\Omega} \langle \nabla E(x-y), \nabla f(x)\rangle\,dx, \quad
\forall y\in {\mathbb R}^N\setminus\overline{\Omega},\ \
\forall p\in [1,\infty]. \end{equation}
For $y\in {\mathbb R}^N\setminus\overline{\Omega}$ fixed,
we define the vector field $Z:\overline{\Omega}\to {\mathbb R}^N$ by
\[Z(x)=f(x)\nabla E(x-y)=\frac{f(x)}{\omega_N |x-y|^N}\,(x-y),\quad
\forall x\in \overline{\Omega}.\]
Clearly, $Z\in C^1(\Omega\setminus A)\cap C(\overline{\Omega})$. Let
$\epsilon>0$ be fixed such that $\overline{B}_\epsilon(a_i)\subset \Omega$,
$\forall i\in I$ and $\overline{B}_\epsilon(a_i)\cap
\overline{B}_\epsilon(a_j)=\emptyset$, $\forall i,j\in I$ with
$i\not=j$. We denote $\Omega_\epsilon:=\Omega\setminus \left(
\cup_{i\in I}\overline{B}_\epsilon(a_i)\right)$. By applying
Proposition~\ref{div} for $Z:\Omega_\epsilon\to {\mathbb R}^N$, we obtain 
\begin{equation} \label{sup}
\begin{aligned}
\int_{\Omega_\epsilon} {\rm div}\,Z(x)\,dx =&\int_{\partial \Omega}
f(x)\frac{\partial E}{\partial\nu}(x-y)\,d\sigma(x)\\
& -\frac{1}{\omega_N}\sum_{i\in I}
\int_{\partial B_\epsilon(a_i)}
\frac{f(x)}{\epsilon |x-y|^{N}}\langle x-y,x-a_i \rangle \,d\sigma(x).
\end{aligned} \end{equation}
Since $y\not\in \overline{\Omega}$, for each $i\in I$, there exists
a constant $M_i>0$ such that
\[ |x-y|>M_i,\quad \forall x\in \partial B_j(a_i),\ \ \forall j\in
(0,\epsilon].\]
Hence, for each $i\in I$, we have
\begin{equation} \label{lab}
\begin{aligned}
& \left|\int_{\partial B_\epsilon(a_i)}
\frac{f(x)}{\epsilon |x-y|^{N}}\langle x-y,x-a_i \rangle \,d\sigma(x)
\right| \leq \int_{\partial B_\epsilon(a_i)} \frac{|f(x)|}{|x-y|^{N-1}}
\,d\sigma(x)\\
& \qquad \leq \frac{\|f\|_{L^\infty(\Omega)}}{M_i^{N_1}}\,{\rm meas}\,(\partial
B_\epsilon(a_i))=\frac{\omega_N \|f\|_{L^\infty(\Omega)}}{M_i^{N-1}}\,
\epsilon^{N-1}\to 0\ \ \mbox{as}\ \epsilon \to 0.
\end{aligned}
\end{equation}
By (\ref{sup}) and (\ref{lab}), it follows that
\begin{equation} \label{etc0}
 \lim_{\epsilon \to 0}\int_{\Omega_{\epsilon}} {\rm div}\,Z(x)\,dx=
\int_{\partial\Omega}
f(x)\frac{\partial E}{\partial\nu}(x-y)\,d\sigma(x).
\end{equation}
Since $x\longmapsto E(x-y) $ is harmonic on ${\mathbb R}^N\setminus \{y\}$,
we find that
\begin{equation} \label{etc1}
\begin{split}
{\rm div}\, Z(x)&=\langle \nabla f(x), \nabla E(x-y)\rangle+
f(x)\Delta_x E(x-y)\\
& =\langle \nabla f(x), \nabla E(x-y)\rangle,\quad
\forall x\in \Omega_\epsilon.
\end{split} \end{equation}
We define $\Psi(x)=|x-y|^{1-N}$, for each $x\in \Omega$. Since
$y\not\in \overline{\Omega}$, we have $\Psi\in C(\overline{\Omega})$ so that
$\Psi\in L^m(\Omega)$, $\forall m\in [1,\infty]$.
By H\"{o}lder's inequality, we infer that
\begin{equation} \label{etc2}
\int_\Omega |\langle \nabla f(x),\nabla E(x-y) \rangle|\,dx  \leq
\frac{1}{\omega_N}\|\nabla f\|_{L^p(\Omega)} \|\Psi\|_{L^{p'}(\Omega)}<\infty,
\quad \forall p\in [1,\infty].
\end{equation}
From (\ref{etc0})--(\ref{etc2}), we conclude (\ref{on}).
\end{proof}

\begin{prop} \label{wa}
If $\Omega$ is an open bounded set with $C^1$ boundary and
$f\in C^2(\Omega)\cap C^1(\overline{\Omega})$ such that
$\Delta f\in C(\overline{\Omega})$, then
\begin{equation}
\begin{split}
\int_{\Omega}\langle \nabla E(x-y),\nabla f(x)\rangle\,dx &=
\int_{\partial\Omega}\frac{\partial f}{\partial \nu}(x)E(x-y)\,d\sigma(x)\\
& \quad -
\int_{\Omega}\Delta f(x) E(x-y)\,dx,\ \ \forall y\in {\mathbb R}^N\setminus \partial\Omega.
\label{grr}
\end{split} \end{equation}
\end{prop}

\begin{proof}
If $y\in {\mathbb R}^N\setminus \overline{\Omega}$, then (\ref{grr})
follows by Proposition~\ref{green} (since $x\longmapsto E(x-y)$ belongs to
$C^2(\Omega)\cap C^1(\overline{\Omega})$).

For $y\in \Omega$ fixed, we choose $\epsilon>0$ such that
${\overline B}_{\epsilon}(y)\subset \Omega$. By Proposition~\ref{green}
(applied on $\Omega\setminus B_\epsilon(y)$), we find
\begin{equation} \label{for}
\begin{aligned}
& \int_{\Omega\setminus B_\epsilon(y)}\Delta f(x)E(x-y)\,dx=
\int_{\partial\Omega}\frac{\partial f}{\partial \nu}(x)E(x-y)\,d\sigma(x)\\
& \quad -\int_{\partial B_\epsilon(y)}
\frac{\partial f}{\partial\nu}(x)E(x-y)\,d\sigma(x)
-\int_{\Omega\setminus B_\epsilon(y)}\langle
\nabla f(x),\nabla E(x-y) \rangle \,dx.
\end{aligned}
\end{equation}
Since $x\longmapsto \Delta f(x)E(x-y)$ is integrable on $\Omega$, we have
\begin{equation} \label{for1}
\int_{\Omega}
\Delta f(x)E(x-y)\,dx=\lim_{\epsilon \to 0}
\int_{\Omega\setminus B_\epsilon(y)}\Delta f(x)E(x-y)\,dx.
\end{equation}
On the other hand, using $f\in C^1(\overline{\Omega})$, we deduce (as
in the proof of Theorem~\ref{teo1}) that $x\longmapsto \langle
\nabla f(x),\nabla E(x-y)\rangle $ is integrable on $\Omega$. It follows that
\begin{equation} \label{for2}
\int_{\Omega} \langle \nabla f(x),\nabla E(x-y)\rangle \,dx=
\lim_{\epsilon\to 0}\int_{\Omega\setminus B_\epsilon(y)}
\langle \nabla f(x),\nabla E(x-y)\rangle \,dx.
\end{equation}
Our next step is to prove that
\begin{equation} \label{for3}
\lim_{\epsilon\to 0}
\int_{\partial B_\epsilon(y)}
\frac{\partial f}{\partial\nu}(x)E(x-y)\,d\sigma(x)
=0.
\end{equation}
Indeed, if $N=2$, then we have
\[ \begin{aligned}
\left|
\int_{\partial B_\epsilon(y)}
\frac{\partial f}{\partial\nu}(x)E(x-y)\,d\sigma(x)
\right| &\leq \int_{\partial B_\epsilon(y)}
\left| \frac{\partial f}{\partial\nu}(x)\right| \frac{1}{2\pi}
\left| \ln |x-y|\right|\,d\sigma(x)\\
& \leq -C\epsilon \log \epsilon\to 0\ \ \mbox{as}\ \epsilon\to 0
\end{aligned}\]
resp., if $N>2$ then
\[ \begin{aligned}
\left|
\int_{\partial B_\epsilon(y)}
\frac{\partial f}{\partial\nu}(x)E(x-y)\,d\sigma(x)
\right| &\leq \int_{\partial B_\epsilon(y)}
\left| \frac{\partial f}{\partial\nu}(x)\right| \frac{1}{\omega_N(N-2)|x-y|^{N-2}}
\,d\sigma(x)\\
& \leq C\frac{{\rm meas}\,(\partial B_\epsilon (y))}
{\epsilon^{N-2}}=C\omega_N \epsilon\to 0\ \ \mbox{as}\ \epsilon\to 0
\end{aligned}\]
where, in both cases, $C$ denotes a positive constant. 

Passing to the limit $\epsilon\to 0$ in (\ref{for}) and using
(\ref{for1})--(\ref{for3}), we obtain (\ref{grr}).
\end{proof}

\begin{rem} Under the assumptions of
Proposition~\ref{wa}, Corollary~\ref{c2}
leads to the Green--Riemann representation formula (see [4, \S{2.4}])
\[ \begin{aligned}
f(y)=&\int_{\partial\Omega} f(x)\frac{\partial E}{\partial \nu}(x-y)\,
d\sigma(x)-\int_{\partial\Omega}\frac{\partial f}{\partial \nu}(x)
E(x-y)\,d\sigma(x)\\
& +\int_{\Omega} \Delta f(x)E(x-y)\,dx,\qquad \forall y\in \Omega
\end{aligned}
\]
and
\[ \begin{aligned}
0=&\int_{\partial\Omega} f(x)\frac{\partial E}{\partial \nu}(x-y)\,
d\sigma(x)-\int_{\partial\Omega}\frac{\partial f}{\partial \nu}(x)
E(x-y)\,d\sigma(x)\\
& +\int_{\Omega} \Delta f(x)E(x-y)\,dx,\qquad \forall y\in {\mathbb R}^N
\setminus \overline{\Omega}.
\end{aligned} \]
Moreover, if $\partial\Omega$ is smooth enough (at least $C^2$), then
\[ f(y)=2\int_{\partial\Omega}f(x)\frac{\partial E}{\partial\nu}(x-y)\,d\sigma(x)
-2\lim_{\Omega\ni t\to y}\int_\Omega \langle \nabla E(x-t),
\nabla f(x)\rangle\,dx,\quad \forall y\in \partial \Omega.\]
\end{rem}

\end{document}